\documentclass[12pt]{amsart}

\usepackage[dvips]{epsfig}

\author{Piotr \'Sniady}
\address{Institute of Mathematics,
University of Wroclaw, pl.\ Grunwaldzki 2/4, 50-384 Wroclaw, Poland}
\email{Piotr.Sniady@math.uni.wroc.pl}

\title{Random regularization of Brown spectral measure}

\theoremstyle{plain}
\newtheorem{lemma}{Lemma}
\newtheorem{theorem}[lemma]{Theorem}
\newtheorem{proposition}[lemma]{Proposition}

\theoremstyle{definition}

\theoremstyle{remark}

\newtheorem*{example}{Example}

\newcommand{\A}{{\mathcal{A}}}
\newcommand{\B}{{\mathcal{B}}}

\newcommand{\E}{{\mathbb{E}}}
\newcommand{\El}{{\mathcal{L}}}

\newcommand{\C}{{\mathbb{C}}}
\newcommand{\R}{{\mathbb{R}}}

\newcommand{\M}{{\mathcal{M}}}

\newcommand{\N}{{\mathbb{N}}}

\newcommand{\gwia}{^{\star}}
\newcommand{\random}{\M_N\big(\El^{\infty-}(\Omega) \big)}

\DeclareMathOperator{\tr}{tr}
\DeclareMathOperator{\Tr}{Tr}
\DeclareMathOperator{\Alg}{Alg}

\begin{document}

\begin{abstract}
We generalize a recent result of Haagerup;
namely we show that a convolution with a standard 
Gaussian random matrix regularizes the behavior of 
Fuglede--Kadison determinant and Brown spectral distribution measure.
In this way it is possible to establish a connection between the limit
eigenvalues distributions of a wide class of random matrices and the
Brown measure of the corresponding limits.
\end{abstract}

\maketitle

\section{Introduction}
The problem of determining the joint distribution of the eigenvalues of a 
random matrix $A^{(N)}\in\random$
with a given distribution of entries is usually very hard and has explicit 
solutions only for very limited cases. A partial solution of this problem is 
to consider the limit $N\rightarrow\infty$ and to hope that it is an easier 
problem than the original one.

Questions of this type can often be solved by Voiculescu's theory of
free probability \cite{VDN}: if a sequence $(A^{(N)})$ of random matrices converges in
$\star$--moments to some $x$, where $x$ is an element of a non--commutative
probability space $(\A,\phi)$, then some properties of matrices $A^{(N)}$
(e.g.\ independence of some entries) imply useful properties of the limit object $x$
(e.g.\ some kind of freeness).
If random matrices $A^{(N)}$ are normal then the empirical distribution of
their eigenvalues converges to the spectral measure of $x$. This approach 
turned out to be successful in determining the limit eigenvalues
distribution in many cases (see e.g.\ \cite{Shl}).

The situation is much more complicated if random matrices $A^{(N)}$ and the 
limit object $x$ are not normal. In this case the spectral measure of $x$ 
has to be replaced by a more complicated object, namely by the Brown
measure of $x$. The problem of determining the Brown measure of a given 
operator is still difficult, but it is much easier than the original 
question about the eigenvalues of a random matrix \cite{Lar, HL, BL}. 
However, since the Brown measure does not behave in a continuous way with 
respect to the topology given by $\star$--moments, the distribution of 
eigenvalues of $A^{(N)}$ does not always converge to the Brown measure of $x$.

Surprisingly, in many known cases when we consider a ``reasonable'' or ``generic''
sequence of random matrices the sequence of their Brown measures converges to the
Brown measure of the limit (cf \cite{BL}). In these cases, however, the 
convergence was proved in this way, that the distribution of eigenvalues of 
$A^{(N)}$ was calculated by ad hoc methods and nearly by an accident it 
turned out to converge to the Brown measure of $x$. Therefore one of the 
most interesting problems in the theory of random matrices is to relate the 
asymptotic distribution of eigenvalues of a sequence of random matrices with 
the Brown measure of the limit.

In this article we show that it is possible to add a small random 
correction to a sequence of random matrices $(A^{(N)})$ which converges 
in $\star$--moments almost surely to some element $x\in(\A,\phi)$ in such a way that the
new corrected sequence still converges to $x$ and that furthermore
the empirical eigenvalues distributions of the new sequence converge to the Brown measure of $x$
almost surely. A slightly different version of this result was proved
recently by Haagerup \cite{Haagerup2001}
and plays a key role in his proof of existence of invariant subspaces 
for a large class of operators. The random correction used by him is a
matrix Cauchy distribution, the first moment of which is unbounded, what
makes it unpleasant for applications. For this reason Haagerup's estimates of
the correction were in the $\mathcal{L}^p$ norm with $0<p<1$.
The random correction considered in this article has a nicer form of a
Gaussian random matrix and for this reason we are able to find better
estimates for the correction, namely in the operator norm.

The Gaussianity of the correction allows us also to find the limit empirical
eigenvalues distribution of a wide class of random matrices, which include 
both the well--known examples of the matrix $G^{(N)}$ with suitably 
normalized independent Gaussian entries (the limit eigenvalues distribution 
was computed by Ginibre \cite{Gin} in the sense of density of eigenvalues,
the almost sure convergence of empirical distributions was proved in
unpublished notes of Silverstein and later in more generality by Bai
\cite{Bai97}; the Brown measure of the limit was computed by Haagerup and
Larsen \cite{Lar,HL}), the so--called elliptic ensemble (the limit eigenvalues
distribution was computed by Petz and Hiai \cite{PH} and the Brown measure of
the limit was computed by Haagerup and Larsen \cite{Lar,HL}) and new examples
for which the eigenvalues distribution was not known before and which are of
the form $G^{(N)}+A^{(N)}$, where entries of $G^{(N)}$ and $A^{(N)}$ are
independent (the Brown measure of the limit of such matrices was computed by
Biane and Lehner \cite{BL}).

Results of this article can be also applied \cite{Sni}
to show that $DT$ operators (which were introduced recently by Dykema and Haagerup
\cite{DH}) maximize microstate free entropy \cite{Voi1994} among all
operators having fixed Brown measure and the second moment.

Our method bases on the observation that if a sequence of random matrices $A^{(N)}$
converges in $\star$--moments to $x$ then the Fuglede--Kadison determinants
$\Delta(A^{(N)})$ converge to $\Delta(x)$ as well if we are able to find some
bottom bounds for the smallest singular values of $A^{(N)}$. Since the random correction
considered in this article is given by a certain matrix--valued Brownian motion, hence
we are able to write a system of stochastic differential equations fulfilled by the singular values.
Unfortunately, finding an exact analytic solution to a non--linear stochastic differential
equation is very difficult. We deal with this problem by proving a certain mononicity
property of our equations and hence we are able to find appropriate bottom estimates for the
singular values.

\section{Preliminaries}
\label{sec:definicje}
\subsection{Non--commutative probability spaces}
A non--commutative probability space is a pair $(\A,\phi)$, where $\A$ is a
$C\gwia$--algebra and $\phi$ is a normal, faithful, tracial state on $\A$.
Elements of $\A$ will be referred to as non--commutative random variables and state $\phi$
as expectation value. The distribution of $x\in\A$ is the collection of all
its $\star$--moments $\big(\phi(x^{s_1}\cdots x^{s_n})\big)$, where 
$s_1,\dots,s_n\in\{1,\star\}$.


\subsection{Fuglede--Kadison determinant}
Let a non--commutative probability space $(\A,\phi)$ be given.
For $x\in\A$ we define its Fuglede--Kadison determinant $\Delta(x)$ by (cf 
\cite{FK}) $$\Delta(x)=\exp\left[ \phi( \ln |x|   ) 
\right].$$

\subsection{Brown measure}
Let a non--commutative probability space $(\A,\phi)$ be given.
For $x\in\A$ we define its Brown measure \cite{Brown} to be the 
Schwartz distribution on $\C$ given by
$$\mu_x= \frac{1}{2\pi} \left( \frac{\partial^2}{\partial a^2} +\frac{\partial^2}{\partial
b^2} \right) \ln \Delta[ x-(a+bi)].$$
One can show that in fact $\mu_x$ is a positive probability measure on $\C$.


\begin{example}
The Brown measure of a normal operator has a particularly easy form;
let $x\in\A$ be a normal operator and let $E$ denote its spectral
measure:
$$x=\int_\C z\ dE(z).$$
Then the Brown measure of $x$ is given by 
$$\mu_x(X)=\phi[ E(X)] $$
for every Borel set $X\subseteq \C$ and the following holds:
$$\phi[ x^k (x\gwia)^l] = \int_{\C} z^k \bar{z}^l\ d\mu_x(z) .$$
\end{example}

\subsection{Random matrices}
We have that $(\M_N,\tr_N)$ is a non--commutative probability space, where
$\M_N$ denotes the set of all complex--valued $N\times N$ matrices and
$\tr_N$ (which for simplicity will be also denoted by $\tr$) is the
normalized trace on $\M_N$ given by
$$\tr_N A=\frac{1}{N} \Tr A \qquad \mbox{for } A\in\M_N,$$
and $\Tr$ denotes the standard trace.

The below simple example shows that for finite matrices the Fuglede--Kadison
determinant $\Delta$ and the usual determinant $\det$ are closely related
and gives heuristical arguments that for every Borel set
$X\subset\C$ the Brown measure
$\mu_x(X)$ provides information on the joint ``dimension'' of
``eigenspaces'' corresponding to $\lambda\in X$.
\begin{proposition}
The Fuglede--Kadison determinant of a matrix $A\in \M_N$ with respect to
a normalized trace $\tr$ is given by $$\Delta(A)=\sqrt[N]{| \det A |}.$$

The Brown measure of a matrix $A\in\M_N$ with respect to the state $\tr$
is a probability counting measure
$$\mu_A=\frac{1}{N} \sum_{i=1}^N \delta_{\lambda_i}, $$
where $\lambda_1,\dots,\lambda_N$ are the eigenvalues of $A$ counted with
multiples.
\end{proposition}

In the following we will be interested in studying the random measure
$\omega\mapsto \mu_{A(\omega)}$ for a random matrix $A\in\random$.
This random measure is called the empirical distribution of eigenvalues.




We will use the following convention: we say that random
matrices $A,B \in\random$ are independent if
the family of entries of $A$ and the family of entries of $B$ are independent.


\subsection{Convergence of $\star$--moments}
Let a sequence $A^{(N)}\in\random$ of
random matrices, a non--commutative probability space $(\A,\phi)$
and $x\in\A$ be given.
We say that the sequence $A^{(N)}$ converges to $x$ in $\star$--moments almost surely
if for every $n\in\N$ and $s_1,\dots, s_n\in\{1,\star\}$ we have that
$$\lim_{N\rightarrow\infty} \tr_N \big[ \big(A^{(N)} \big)^{s_1} \cdots
\big(A^{(N)} \big)^{s_n} \big] = \phi( x^{s_1} \cdots x^{s_n} )$$
holds almost surely.

Let a sequence $A^{(N)}\in\random$ of
random matrices, a non--commutative probability space $(\A,\phi)$
and $x\in\A$ be given. We say that the sequence $A^{(N)}$ converges to $x$ in expected
$\star$--moments
if for every $n\in\N$ and $s_1,\dots, s_n\in\{1,\star\}$ we have that
$$\lim_{N\rightarrow\infty} \E \tr_N \big[ \big(A^{(N)} \big)^{s_1} \cdots
\big(A^{(N)} \big)^{s_n} \big] = \phi( x^{s_1} \cdots x^{s_n} ).$$

\subsection{Discontinuity of Fuglede--Kadison determinant and Brown measure}
One of the greatest difficulties connected with the Fuglede--Kadison  
determinant
and Brown spectral distribution measure is that---as we shall see in the 
following example---these 
two objects do not behave in a continuous way with respect to the topology
given by convergence of $\star$--moments.

We say that $u\in\A$ is a Haar unitary if $u$ is unitary and 
$\phi(u^k)=\phi\big( (u\gwia)^k \big)=0$ for every $k=1,2,\dots$
It is not difficult to see that the sequence $(\Xi^{(N)})$ converges in
$\star$--moments to the Haar unitary, where $\Xi^{(N)}$ is an $N\times N$
nilpotent matrix
\begin{equation} \Xi^{(N)}=\left[ \begin{array}{ccccc} 
0       & 0     &\cdots & 0     & 0             \\
1       & 0     &\cdots & 0     & 0             \\
0       & 1     &\ddots &\vdots & \vdots        \\      
\vdots  &\ddots &\ddots & 0     & 0             \\
0       &\cdots & 0     & 1     & 0
\end{array}\right].
\label{eq:nilpotent}
\end{equation}
Every matrix $\Xi^{(N)}$ has the determinant equal to 
$0$, while the Haar unitary has the Fuglede--Kadison determinant equal to 
$1$; every matrix $\Xi^{(N)}$ has the Brown measure equal to $\delta_0$, 
while the Brown measure of the Haar unitary is the uniform measure on the 
unit circle $\{z\in\C:|z|=1\}$.

The reason for the discontinuity of Fuglede--Kadison determinant is that
the logarithm is not bounded from below on any interval $[0,t]$. However, since it
is bounded from above, Fuglede--Kadison determinant is upper--semicontinuous.

\begin{lemma}  \label{lem:kopalny}
Let $A^{(N)}$ be a sequence of random matrices which converges in $\star$--moments to
a non--commutative random variable $x$ almost surely.
Then for every $\lambda\in\C$
$$\limsup_{N\rightarrow\infty} \tr \ln |A^{(N)}-\lambda| \leq \ln \Delta(x-\lambda)$$
holds almost surely.

Let $A^{(N)}$ be a sequence of random matrices which converges in expected $\star$--moments
to a non--commutative random variable $x$.
Then for every $\lambda\in\C$ we have
$$\limsup_{N\rightarrow\infty} \E \tr \ln |A^{(N)}-\lambda| \leq \ln \Delta(x-\lambda).$$
\end{lemma}
\begin{proof}
For each $\epsilon>0$ there exists an
even polynomial $Q$ such that
$$\ln r\leq Q(r) \qquad \text{for every } r>0$$
and
$$Q(r) \leq \frac{\ln (r^2+\epsilon)}{2} \qquad \text{for every } 0\leq r \leq \|x\|.$$
Hence
$$\tr \ln |A^{(N)}-\lambda| \leq \tr Q\big( |A^{(N)}-\lambda| \big).$$
The right--hand side converges almost surely (resp.\ in the expectation value) to
$\phi\big[ Q\big( |x-\lambda| \big) \big]\leq \phi\left(  \frac{\ln (r^2+\epsilon)}{2}  \right)$.
By taking the limit $\epsilon\rightarrow 0$ both parts of the lemma follow.
\end{proof}

\subsection{Gaussian random matrices}
We say that a random matrix 
$$G^{(N)}=(G^{(N)}_{ij})_{1\leq i,j\leq N}\in\random$$
is a standard Gaussian random matrix if 
$$\big(\Re G^{(N)}_{ij}\big)_{1\leq i,j\leq N}, 
\big(\Im G^{(N)}_{ij}\big)_{1\leq i,j \leq N}$$ 
are independent Gaussian variables with mean zero and variance 
$\frac{1}{2 N}$. 

We say that $$M^{(N)}:\R_{+}\rightarrow \random, \qquad
M^{(N)}(t)=\big(M^{(N)}_{ij}(t) \big)_{1\leq i,j\leq N}$$ 
is a standard matrix Brownian 
motion if 
$$\big( \Re M^{(N)}_{ij} \big)_{1\leq i,j\leq N}, \big( \Im  M^{(N)}_{ij}\big)_{1\leq i,j 
\leq N}$$ 
are independent Brownian motions which are normalized in such a 
way that the variance is given by
$$\E \big( \Re M^{(N)}_{ij}(t) \big)^2 = \E \big( \Im M^{(N)}_{ij}(t) 
\big)^2= \frac{t}{2 N}.$$


\subsection{Circular element}
There are many concrete characterizations of the Voiculescu's circular
element $c$ \cite{VDN} but we will use the following implicit definition.
One can show that the sequence $G^{(N)}$ converges both in expected $\star$--moments
and in $\star$--moments almost surely to a certain
non--commutative random variable $c$ \cite{Voi1991,Tho}.

\subsection{Freeness}
Let $(\A,\phi)$ be a non--commutative probability space and let
$(\A_i)_{i\in I}$ be a family of unital $\star$--subalgebras of $\A$. We say
that the algebras $(\A_i)_{i\in I}$ are free if
$$\phi(x_1 x_2 \cdots x_n)=0 $$
holds for every $n\geq 1$, every $i_1,i_2,\dots,i_n\in I$ such that $i_1\neq
i_2$, $i_2\neq i_3$,\dots, $i_{n-1}\neq i_n$, and every
$x_1\in\A_{i_1},\dots, x_n\in\A_{i_n}$ such that
$\phi(x_1)=\cdots=\phi(x_n)=0$ (cf \cite{VDN}).

Let $(X_i)_{i\in I}$ be a family of subsets of $\A$. We say that sets $X_i$
are free if unital $\star$--algebras $(\Alg \{X_i, X_i\gwia\})_{i\in I}$ are
free.


\section{Singular values of matrix Brownian motions}
%
Let $A\in\random$ be a given random
matrix and $M^{(N)}$ be a standard matrix Brownian motion
such that $A$ and $M^{(N)}$ are independent.  For $t\geq 0$ we define
a random matrix $A_{t}$ by
$$ A_{t}(\omega)=A(\omega)+M^{(N)}(t,\omega). $$
It should be understood that matrix Brownian motions $(M^{(N)})_{N=1,2,\dots}$ are independent.

If we are interested in $A_t$ for only one value of $t\geq 0$ we can
express $A_t$ as follows:
$$A_{t}(\omega)=A(\omega)+\sqrt{t}\ G^{(N)}(\omega),$$
where $G^{(N)}$ is a standard Gaussian random matrix such that $A$
and $G^{(N)}$ are independent.

If $x\in\A$ is a non--commutative random variable, we
we can always extend the algebra $\A$ and find
$c\in\A$ such that $\{x,x\gwia\}$ and $\{c,c\gwia\}$ are free and $c$
is a circular element \cite{VDN}. We will denote
$$x_t= x+ \sqrt{t}\ c.$$

\begin{proposition}
\label{prop:steen}
If sequence of random matrices $|A^{(N)}|^2$ converges in $\star$--moments to
$|x|^2$ almost surely then for every $t\geq 0$
the sequence $|A^{(N)}_t|^2$ converges in $\star$--moments to $|x_t|^2$ almost surely.

If sequence of non--random matrices $|A^{(N)}|^2$ converges in $\star$--moments to
$|x|^2$ then for every $t\geq 0$
the sequence $|A^{(N)}_t|^2$ converges in expected $\star$--moments to $|x_t|^2$.
\end{proposition}
\begin{proof}
The first part of the propositions follows under additional assumption that $\sup_N \|A^{(N)}\|<\infty$
almost surely from recent results of Hiai and Petz \cite{HP}. For the general case
observe that since for all unitary matrices $U,V\in\M_N$ and $n\in\N$ the distributions of random variables
$\tr |A^{(N)}+\sqrt{t} G^{(N)}|^{2n}$ and $\tr |A^{(N)}+\sqrt{t} U G^{(N)} V|^{2n}=
tr |V\gwia A^{(N)} U\gwia+\sqrt{t} G^{(N)}|^{2n}$ coincide,
hence it is enough to prove the first part under assumption that every matrix $A^{(N)}$ is almost
surely diagonal. The method of Thorbj\o{}rnsen can be generalized to this case \cite{Tho}.

The second part of the proposition was proved by Voiculescu \cite{Voi1991}.
\end{proof}

For any
$t\geq 0$ and $\omega\in\Omega$ let $\lambda_1(t,\omega)\geq \cdots \geq
\lambda_N(t,\omega)$ denote singular values of the matrix
$A_t(\omega)$.

In Section \ref{subsec:wyprowadzenie} we derive
stochastic differential equations for $\lambda_1,\dots,\lambda_N$
using similar methods to those of Chan \cite{Chan} and obtain
\begin{equation}
d\lambda_i(t)=\Re (dB_{ii})+
\frac{dt}{2 \lambda_i} \left(1-\frac{1}{2 N}+
\sum_{j\neq i} \frac{\lambda_i^2+\lambda_j^2}{N (\lambda_i^2-\lambda_j^2)}
\right),
\label{eq:rownanieDysonaPrim}
\end{equation}
where $B$ is a standard matrix Brownian motion.

\begin{theorem}
\label{theo:twierdzenieoporownywaniu}
Let $A^{(1)}$ and $A^{(2)}$ be non--random matrices of the same size,
$A^{(1)}, A^{(2)}\in\M_N$. For $n=1,2$ let $s^{(n)}_1 \geq
\cdots \geq s^{(n)}_N$ be the singular values of the matrix $A^{(n)}$.
Suppose that for each $1\leq k\leq N$ we have $s^{(1)}_k< s^{(2)}_k$.

Then for every $t\geq 0$ there exists a probability space 
$(\Omega, {\mathcal{B}},P)$ and random matrices $G^{(1)}, G^{(2)}\in\random$
such that each matrix $G^{(i)}$ is
a standard Gaussian random matrix (but matrices $G^{(1)}$ and $G^{(2)}$
might be dependent) and such that
$$\tr f\Big( \big| A^{(1)}+\sqrt{t}\ G^{(1)}(\omega) \big| \Big) \leq
\tr f\Big( \big| A^{(2)}+\sqrt{t}\ G^{(2)}(\omega) \big| \Big)$$
holds for every $\omega\in\Omega$ and every nondecreasing function
$f:\R\rightarrow\R$. 
\end{theorem}
\begin{proof}
Let us consider a probability space $(\Omega,\B,P)$, a standard matrix
Brownian motion $B:\R_+\rightarrow \random$
and for each $n\in\{1,2\}$ we find the solution of the system of stochastic 
differential equations
\begin{equation}
d\lambda^{(n)}_i=  
\Re (dB_{ii})+
\frac{dt}{2 \lambda^{(n)}_i} \left(1-\frac{1}{2 N}+
\sum_{j\neq i} \frac{(\lambda^{(n)}_i)^2+(\lambda_j^{(n)})^2}{N 
((\lambda_i^{(n)})^2-(\lambda_j^{(n)})^2)} \right).
\label{eq:rownanieDysonaPrim2} \end{equation}
together with initial conditions 
$$\lambda_i^{(n)}(0,\omega)=s_i^{(n)}.$$
We have that for each $t>0$ and $n\in\{1,2\}$ the joint distribution of random variables
$\lambda_i^{(n)}(t)$, $i=1,\dots,N$ coincides with the joint distribution
of singular values of the random matrix $A^{(n)}_t$.

The theorem will follow from the following stronger statement:
for almost every $\omega$ and every $t\geq 0$ we have
\begin{equation}
\lambda^{(1)}_i(t,\omega)< \lambda^{(2)}_i
(t,\omega) \qquad \mbox{for every } 1\leq i\leq N.
\label{eq:warunek}
\end{equation}

{} From Eq.\ (\ref{eq:rownanieDysonaPrim2}) it follows that for almost every 
$\omega\in\Omega$ we have that $\lambda^{(1)}-\lambda^{(2)}$ has a continuous
derivative (see Section \ref{subsubsec:avoid}).
For a fixed $\omega\in\Omega$ let $t_0$ be the smallest $t\geq 0$ such that
(\ref{eq:warunek}) does not hold. Trivially we have $t_0>0$. There exists
an index $j$ such that $\lambda^{(1)}_j(t_0)=\lambda^{(2)}_j(t_0)=:\lambda_j$
and for every $i$ we have $\lambda^{(1)}_i(t_0)\leq
\lambda^{(2)}_i(t_0)$.  Eq.\ (\ref{eq:rownanieDysonaPrim2}) gives us
\begin{equation}
\frac{d}{dt} \Big(\lambda_j^{(1)}(t)-\lambda_j^{(2)}(t) \Big)\bigg|_{t=t_0}=
\sum_{k\neq j} \frac{\lambda_j \Big( (\lambda_k^{(1)})^2 - (\lambda_k^{(2)})^2 \Big)}{
 N \Big( \lambda_j^2 - (\lambda_k^{(1)})^2 \Big) \Big( \lambda_j^2-(\lambda_k^{(2)})^2 \Big)}.
\label{eq:rownanieDysonaPrim3}
\end{equation}

It is easy to see that if there exists at least one index $1\leq
k\leq N$ such that $\lambda^{(1)}_k(t_0)\neq \lambda^{(2)}_k(t_0)$ then
$$\frac{d}{dt}\Big( \lambda_j^{(1)}-\lambda_j^{(2)} \Big)\bigg|_{t=t_0}<0,$$
so it follows that for small $d>0$ we have 
$\lambda_j^{(1)}(t)-\lambda_j^{(2)}(t)>0$ for $t_0-d<t<t_0$. This contradicts the
minimality of $t_0$.

We define $\delta(t)=\lambda_i^{(1)}(t)-\lambda_i^{(2)}(t)$. If we replace in
(\ref{eq:rownanieDysonaPrim3}) $\lambda^{(2)}_i$ by $\lambda^{(1)}_i-\delta_i$
then it becomes a system of non--stochastic ordinary differential equations for
$\delta_i$. If for all indexes $1\leq i\leq N$ we have
$\lambda^{(1)}_i(t_0)=\lambda^{(2)}_i(t_0)$ then $\delta_i(t_0)=0$ and the solution
exists and is unique in some (backward) interval. This contradicts the minimality of $t_0$.
\end{proof}

\begin{proposition}
\label{prop:niespada}
If $A^{(N)}\in\random$ is a random
 matrix and $\lambda\in\C$ then the function
$\R_+ \ni t\mapsto\E \tr \ln |A^{(N)}_t-\lambda|$ is nondecreasing.

For $x\in\A$ and $\lambda\in\C$
we have that the function $\R_+\ni t\mapsto \ln \Delta(x_t)$ is nondecreasing and
\begin{equation} \label{eq:dokrajutegogdziekruszynechleba}
\lim_{t\rightarrow 0^{+}} \ln \Delta(x_t) = \ln \Delta(x).
\end{equation}
\end{proposition}
\begin{proof}
We can regard $\M_N$ as a $2N^2$--dimensional real
Euclidean space equipped with a scalar product
$\langle m,n \rangle=\Re \Tr m n\gwia$.
As usually we define the Laplacian to be $\nabla^2=\sum_{1\leq k\leq
2N^2} D_{v_k}^2$, where $v_1,\dots,v_{2N^2}$ is the
orthonormal basis of this space and $D_v$ is a derivative
operator in direction $v$.

Notice that $\ln | \det A |=\Re \ln \det A$. We can regard $\det A$ as
a holomorphic function of $N^2$ complex variables (=entries of the matrix).
On the other hand it is a known--fact that if $f(z_1,\dots,z_k)$ is a
holomorphic function then the Laplacian of its logarithm
is a positive measure. This and It\^o formula imply the first part of the proposition.

For the second part we construct a sequence $(A^{(N)})$, where $A^{(N)}\in\M_N$,
such that $A^{(N)}$ converges in $\star$--moments to
$|x|$ and such that $\lim_{N\rightarrow\infty} \tr \ln (A^{(N)})=
\ln \Delta(x)$ and apply Lemma \ref{lem:kopalny} for the sequence $A^{(N)}_t$.
This shows that $\R_+\ni t\mapsto \ln \Delta(x_t)$ is nondecreasing.

On the other hand the inequality
$$\limsup_{t \rightarrow 0^+} \ln \Delta(x_t) \leq \ln \Delta(x)$$
can be proved similarly as in Lemma \ref{lem:kopalny}.
\end{proof}


\section{The main result}
\begin{theorem}
\label{theo:glowne}
Let $A^{(N)}\in\random$ be a
sequence of random  matrices such that $A^{(N)}$ converges in $\star$--moments
to a non--commutative random variable $x$ almost surely.

For every $t>0$ we have that the sequence of empirical distributions $\mu_{A^{(N)}_t(\omega)}$ converges
in the weak topology to $\mu_{x_t}$ almost surely.

There exists a sequence $(t_N)$ of positive numbers such that
$\lim_{N\rightarrow\infty} t_N=0$ and
the sequence of empirical distributions $\mu_{A^{(N)}_{t_N}(\omega)}$ converges
in the weak topology to $\mu_{x}$ almost surely.
\end{theorem}

\begin{theorem}
\label{theo:glowne2}
Let $(A^{(N)})$ be a sequence of non--random matrices ($A^{(N)}\in\M_N$)
which converges in $\star$--moments to a non--commutative random
variable $x\in\A$, where $(\A,\phi)$ is a non--commutative probability space.

There exists a sequence $(\tilde{A}^{(N)})$ of non--random matrices
such that the distributions of eigenvalues $\mu_{\tilde{A}^{(N)}}$ converge
weakly to $\mu_x$ and
$$\lim_{N\rightarrow\infty} \| A^{(N)}-\tilde{A}^{(N)} \| =0, $$
where $\|\cdot\|$ denotes the operator norm of a matrix.
\end{theorem}

\begin{figure}  
\centering
\psfig{file=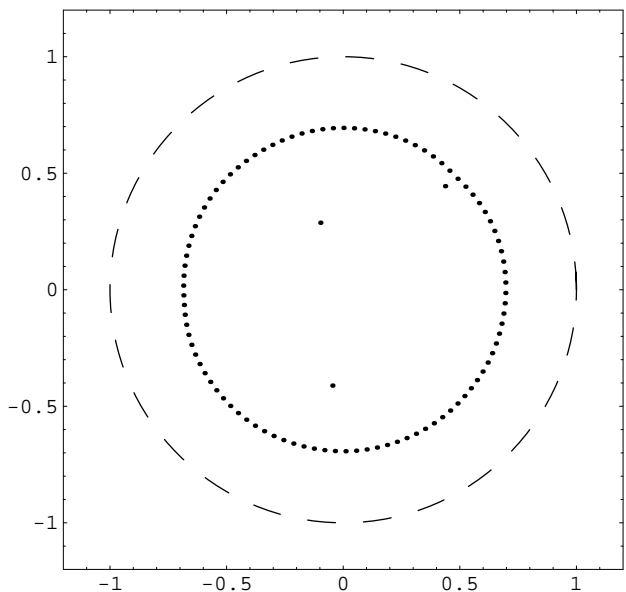,width=0.65 \textwidth}
\caption{Sample eigenvalues of a random matrix $\Xi^{(N)}+\sqrt{t} G^{(N)}$ for
$N=100$ and $t=10^{-100}$.}
\label{fig:fig1}
\end{figure}

\begin{figure}  
\centering
\psfig{file=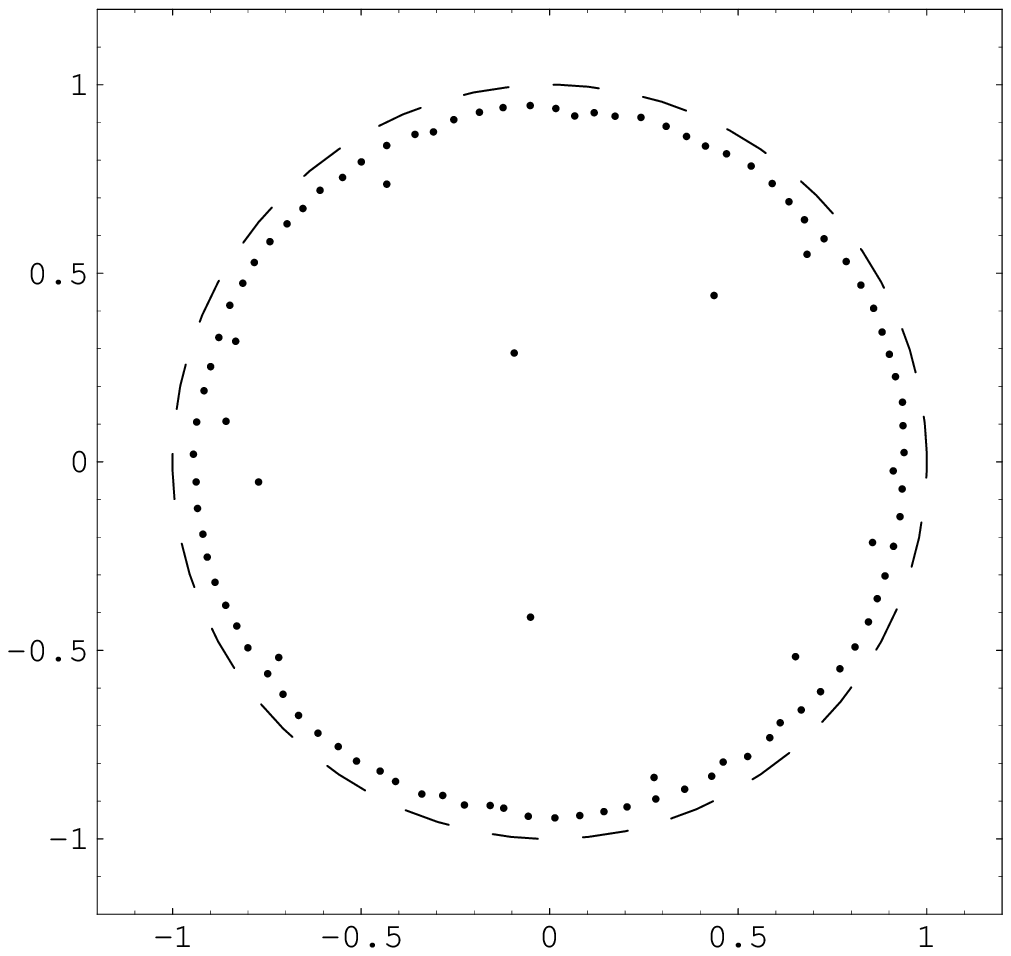,width=0.65 \textwidth}
\caption{Sample eigenvalues of a random matrix $\Xi^{(N)}+\sqrt{t} G^{(N)}$ for
$N=100$ and $t=10^{-5}$.}
\label{fig:fig2}
\end{figure}

\begin{figure}  
\centering
\psfig{file=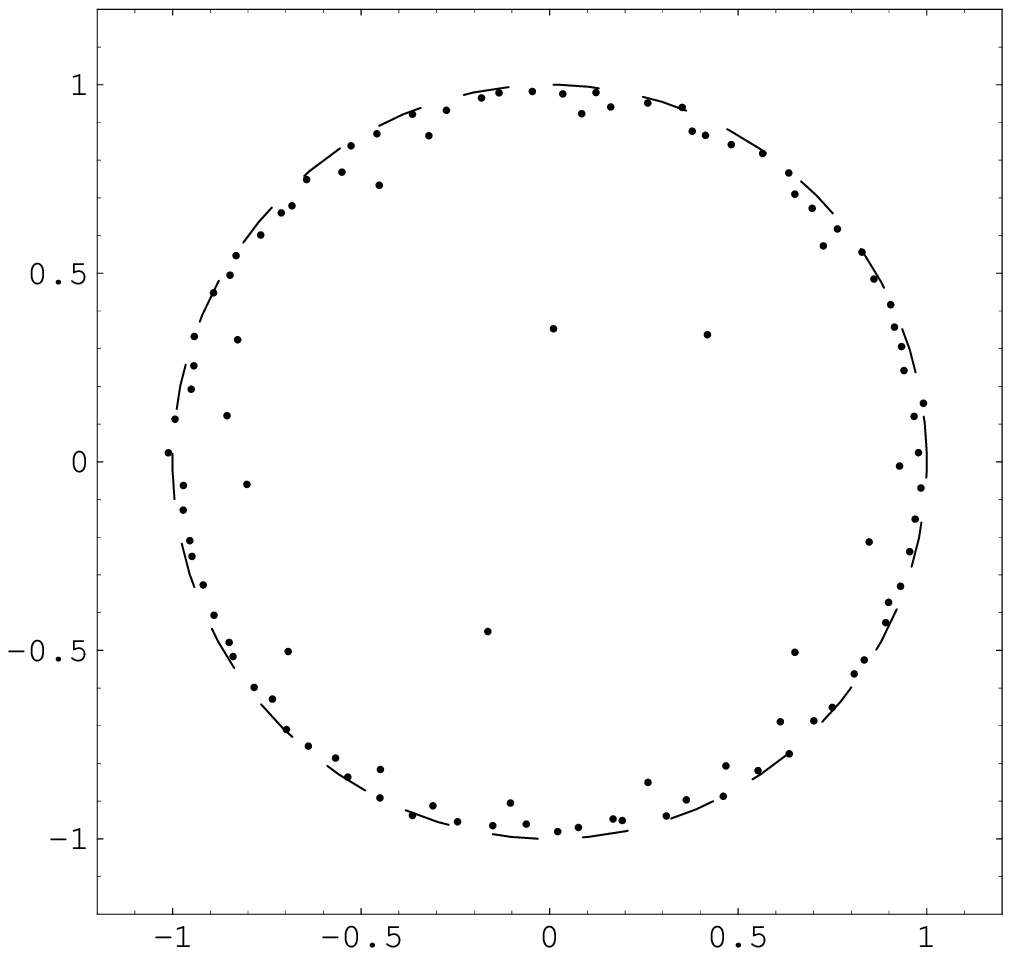,width=0.65 \textwidth}
\caption{Sample eigenvalues of a random matrix $\Xi^{(N)}+\sqrt{t} G^{(N)}$ for
$N=100$ and $t=10^{-2}$.}
\label{fig:fig3}
\end{figure}

\begin{figure}  
\centering
\psfig{file=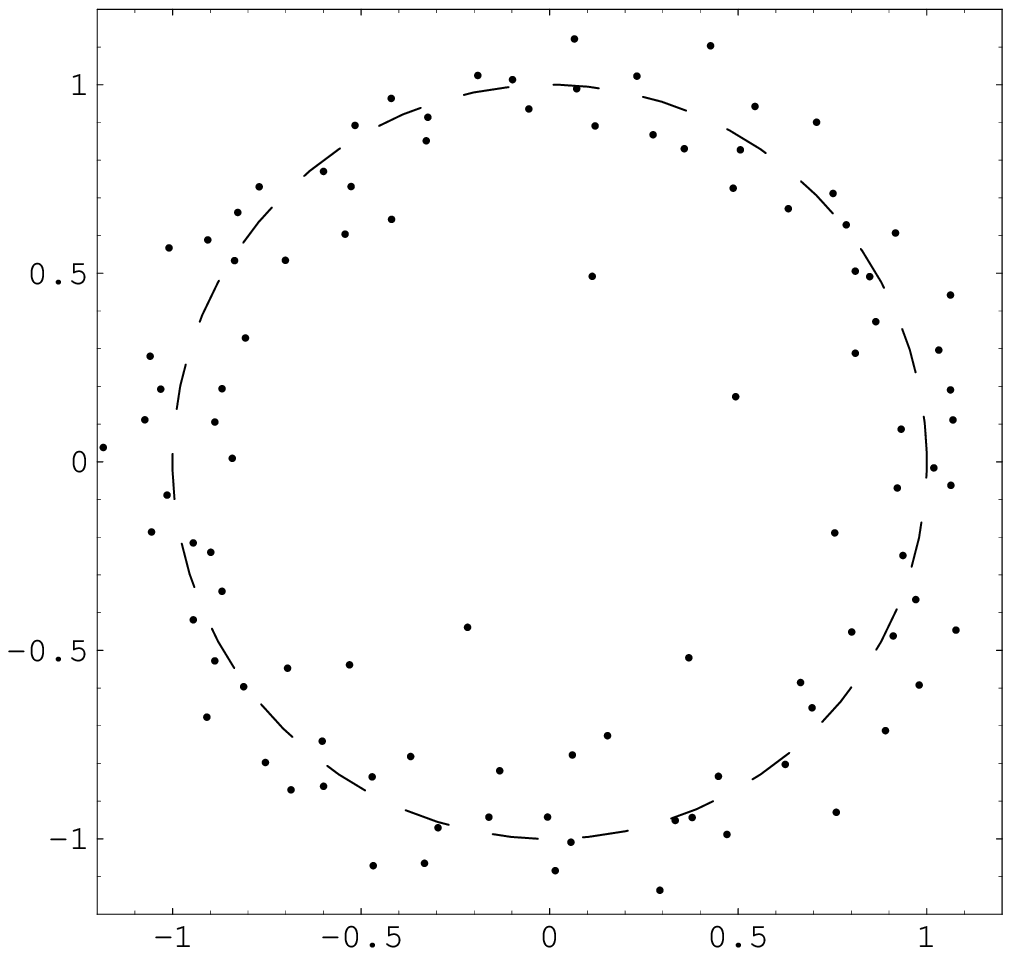,width=0.65 \textwidth}
\caption{Sample eigenvalues of a random matrix $\Xi^{(N)}+\sqrt{t} G^{(N)}$ for
$N=100$ and $t=3 \cdot 10^{-1}$.}
\label{fig:fig4}
\end{figure}

As an illustration to the above theorems
we present on Fig.\ \ref{fig:fig1}---\ref{fig:fig4} results of a computer
experiment; we plotted eigenvalues of the nilpotent matrix $\Xi_N$ from Eq.\
(\ref{eq:nilpotent}) with a random Gaussian correction. The size of the
matrices was $N=100$; with dashed line we marked the spectrum of the Haar
unitary, which is the circle of radius $1$ centered in $0$. We recall that
the sequence $(\Xi^{(N)})$ converges in $\star$--moments to the Haar unitary.
As one can see if the random correction is too small then the eigenvalues
of the corrected matrix behave like the the eigenvalues of $\Xi^{(N)}$
and if the random correction is too big then the eigenvalues of the corrected
matrix are dispersing on the plane.

An interesting problem for future research is for a given sequence $A^{(N)}$ of
random matrices which converges in $\star$--moments to $x$
and fixed $N$ to determine the optimal value of $t_N$ for which
the measure $\mu_{A^{(N)}_{t_N}}$ is the best approximation of $\mu_x$

Before we present the proofs of these theorems we shall prove the following
lemma.

\begin{lemma}  \label{lem:kukuruku}
Let $A^{(N)}$ be as in Theorem \ref{theo:glowne}. For every $t>0$ and
 every $\lambda\in\C$ we have that
$$\lim_{N\rightarrow\infty} \tr \ln |A_t^{(N)}-\lambda| = \ln \Delta(x_t-\lambda)$$
holds almost surely.
\end{lemma}
\begin{proof}
Let us fix $\lambda\in\C$. For any $\epsilon>0$ we define functions on $\R_+$
$$f_\epsilon(r)= \frac{\ln (r^2+\epsilon)}{2}, $$
$$g_\epsilon(r)= \ln r-\frac{\ln (r^2+\epsilon) }{2}.$$
Each function $f_\epsilon$ is well defined on $[0,\infty)$ and
$f_\epsilon$ converges to the function $t\mapsto \ln t$ pointwise as $\epsilon$ tends to $0$.
Each function $g_\epsilon$ is increasing and $g_\epsilon$
converges pointwise to $0$ as $\epsilon$ tends to $0$.

Since the function $f_\epsilon$ has a polynomial growth at infinity,
therefore there exist even polynomials $S(r)$ and $Q(r)$ such that
$S(r)\leq f_\epsilon(r) \leq Q(r)$ holds for every $r\geq 0$
and furthermore $| S(r)- f_\epsilon(r)|<\epsilon$, $|Q(r) - f_\epsilon(r)|<\epsilon$
hold for every $0\leq r\leq \|x_t\|$.

We apply Theorem \ref{theo:twierdzenieoporownywaniu} for a pair of matrices
$0$ and $A^{(N)}-\lambda$ and obtain a probability space $(\Omega,\B,P)$ and
Gaussian random matrices $\tilde{G}^{(N)},G^{(N)}\in \random$ such that
$$\tr g_\epsilon (|\sqrt{t}\ \tilde{G}^{(N)}|) \leq
\tr g_\epsilon(|A^{(N)}_t-\lambda|)$$
holds for every $\omega\in\Omega$, where as 
usually $A_t^{(N)}=A^{(N)}+\sqrt{t}\ G^{(N)}$. For simplicity here and in 
the following we skip the obvious dependence of random variables on $\omega$.

We have that
\begin{multline*}\tr \ln |A_t^{(N)}-\lambda|= \tr f_\epsilon
\big(|A_t^{(N)}-\lambda| \big)+ \tr g_\epsilon\big( |A_t^{(N)}-\lambda|
\big) \geq  \\
\tr S(|A_t^{(N)}-\lambda|)  + \tr g_\epsilon
(|\sqrt{t}\ \tilde{G}^{(N)}|) =\\
\tr S(|A_t^{(N)}-\lambda|)  + \tr \ln (|\sqrt{t}\ \tilde{G}^{(N)}|) -
\tr f_\epsilon (|\sqrt{t}\ \tilde{G}^{(N)}|) \geq \\
\tr S(|A_t^{(N)}-\lambda|)  + \tr \ln (|\sqrt{t}\ \tilde{G}^{(N)}|) -
\tr Q(|\sqrt{t}\ \tilde{G}^{(N)}|)=:X^{(N)}.
\end{multline*}
Propositions \ref{prop:steen} and \ref{prop:wyznacznikgausa} show that
\begin{multline*}
\lim_{N\rightarrow\infty} X^{(N)}= \phi\big( S(|x_t-\lambda|) \big) + \phi( \ln |\sqrt{t} c|)-
\phi\big( Q(|\sqrt{t} c|) \big) \geq \\
\phi(\ln |x_t-\lambda|)+\phi\big( g_\epsilon(|\sqrt{t}\ c|) \big)- 2 \epsilon
\end{multline*}
holds almost surely. Hence by taking the limit $\epsilon\rightarrow\infty$ we obtain that
the inequality
$$\liminf_{N\rightarrow\infty} \tr \ln |A_t^{(N)}-\lambda| \geq \ln \Delta(x_t-\lambda)$$
holds almost surely.

The upper estimate
$$\limsup_{N\rightarrow\infty} \tr \ln |A_t^{(N)}-\lambda| \leq \ln \Delta(x_t-\lambda)$$
follows from Lemma \ref{lem:kopalny} and Proposition \ref{prop:steen}, what finishes the proof.
\end{proof}

\begin{proof}[Proof of Theorem \ref{theo:glowne}]
For the proof of the first part of the theorem let
let $K\subset\C$ be a compact set.  In the following $y$ will denote either
$x\in\A$ or the matrix $A_{t}^{(N)}(\omega)$.
Let $f\in C^2(K)$ be a smooth enough function with a compact support
$K\subset\C$. From the definition of the Brown measure we have
\begin{multline*}\int_{\C} f(\lambda)\ d\mu_y(\lambda)=\frac{1}{2\pi} \langle f(\lambda),
\nabla^2 \ln \Delta(y-\lambda) \rangle =
\\ \frac{1}{2\pi} \int_{\C} \ln \Delta(y-\lambda) \nabla^2
f(\lambda)\ d\lambda
\end{multline*}
Since twice differentiable functions $C^2(K)$ are dense in the
set of all continuous functions $C(K)$ therefore the almost certain convergence of
measures $\mu_{A^{(N)}_{t}(\omega)}$ in the weak topology
to the measure $\mu_{x}$ would follow if the
sequence of functions $\tr \ln \big|A_{t}^{(N)}(\omega)-\lambda\big|$ converges to the function
$\ln \Delta(x-\lambda)$ in the local $\El^1$ norm almost surely. Therefore it would be sufficient
to show that for almost every $\omega\in\Omega$ we have (for simplicity here and in
the following we skip the obvious dependence of random variables on $\omega$)
$$\lim_{N\rightarrow\infty} \int_K \Big| \tr\ln |A_t^{(N)}-\lambda| - \ln \Delta(x_t-\lambda) \Big|
d\lambda=0.$$

From Lemma \ref{lem:kukuruku} and the Fubini theorem
follows that for almost every $\omega\in \Omega$ we have
$$\lim_{N\rightarrow\infty} \tr\ln |A^{(N)}_t-\lambda| = \ln \Delta(x_t-\lambda)$$
for almost all $\lambda\in K$. Now it is sufficient to show that
\begin{equation}
\int_K \sup_N \Big| \tr\ln |A^{(N)}_t-\lambda| \Big| d\lambda+
\int_K \Big|\ln \Delta(x_t-\lambda) \Big| d\lambda < \infty
\label{eq:januszkiewicz}
\end{equation}
holds almost surely in order to apply the majorized convergence theorem.

Note that $\lambda\mapsto\log \Delta(x_t-\lambda)$
is subharmonic \cite{Brown} and hence it is a local ${\mathcal{L}}^1$ function; therefore
we only need to find estimates for the first summand in (\ref{eq:januszkiewicz}).

Theorem \ref{theo:twierdzenieoporownywaniu} gives us that for almost every $\omega\in\Omega$
$$ \tr \ln |A_t^{(N)}-\lambda| \geq \tr \ln |\sqrt{t} \tilde{G}^{(N)}|  $$
hence Proposition \ref{prop:wyznacznikgausa} implies that
$$ \E \min\Big(0, \inf_N \tr \ln |A_t^{(N)}-\lambda|\Big)
 \geq \E \min\Big(0, \inf_N \tr \ln |G^{(N)}|+\frac{\ln t}{2}\Big)$$
is uniformly bounded from below over $\lambda\in\C$.
From Fubini theorem follows that
for almost every $\omega\in\Omega$ we have
$$\int_K  \min\Big(0,\inf_N \tr \ln |A^{(N)}_t(\omega)-\lambda|\Big) d\lambda > -\infty.$$

From the simple inequality $\log r<r^2$ which holds for every $r>0$ we have
$$\tr \ln |A^{(N)}_t-\lambda| < \tr |A^{(N)}_t-\lambda|^2\leq \sqrt{\tr |A_t^{(N)}|^2+\lambda^2}.$$
By Proposition \ref{prop:steen} we have that $\tr |A_t^{(N)}|^2$ converges almost surely,
hence the family of functions $K\ni\lambda\mapsto \tr \ln |A^{(N)}_t-\lambda|$ is almost surely
uniformly bounded from above, what finishes the proof of the first part of the theorem.

From the first part of theorem follows that
there exists a decreasing sequence $(t_N)$ of positive numbers which converges to $0$ and such that
for any compact $K\subset\C$
$$
\lim_{N\rightarrow\infty} \int_K \Big| \tr\ln |A_{t_N}^{(N)}-\lambda| - \ln \Delta(x_{t_N}-\lambda) \Big|
d\lambda=0 $$
holds almost surely.

Proposition \ref{prop:niespada} implies that the majorized convergence
theorem can be applied (we recall that $\lambda\mapsto\log \Delta(y-\lambda)$
is always a local ${\mathcal{L}}^1$ function) hence
$$\lim_{N\rightarrow\infty} \int_K \Big| \ln \Delta(x_{t_N}-\lambda) - \ln \Delta(x-\lambda) \Big| d\lambda=0.$$

The above two equations combine to give
$$\lim_{N\rightarrow\infty} \int_K \Big| \tr\ln |A_{t_N}^{(N)}-\lambda| - \ln \Delta(x-\lambda) \Big|
d\lambda=0 $$
almost surely. The convergence of empirical distributions of eigenvalues follows now exactly as in the
proof of the first part.
\end{proof}

\begin{proof}[Proof of Theorem \ref{theo:glowne2}]
Let $(t_N)$ be a sequence given by Theorem \ref{theo:glowne}. Since
$\limsup_{N\rightarrow\infty} \|G^{(N)}\|<\infty$ almost surely \cite{Gem}, hence for
almost every $\omega\in\Omega$
$$\tilde{A}^{(N)}=A^{(N)}+\sqrt{t_N} G^{(N)}(\omega) $$
is the wanted sequence.
\end{proof}


\section{Technical results}
\subsection{Derivation of the stochastic differential equation for $\lambda_i$}
\label{subsec:wyprowadzenie}
\subsubsection{Singular values as functions on $\M_N$}
In this subsection we are going to evaluate the first and
the second derivative of the map
$$s:\M_N\ni m\mapsto \big( s_1(m),\dots,s_N(m) \big),$$
where $s_1(m),\dots,s_N(m)$ denote the singular values of a matrix $m$.

The perturbation theory shows (cf chapter II.2 of \cite{Kat})
that if $D$ is a diagonal matrix with
eigenvalues $\nu_1,\dots,\nu_N$ such that $\nu_i\neq\nu_j$ for all $i\neq j$
and $\Delta D$ is any matrix then the eigenvalues
$\nu'_1,\dots,\nu'_N$ of a matrix $D+\Delta D$ are given by
$$ \nu'_i=\nu_i + \Delta D_{ii} + \sum_{j\neq i} \frac{\Delta D_{ij}
\Delta D_{ji}}{\nu_i-\nu_j}+ O\big(\| \Delta D\|^3 \big)  $$
for small enough $\| \Delta D\|$ and that the map
$\Delta D\mapsto (\nu'_1,\dots,\nu'_N)$ is $C^2$ in some neighbourhood of $0$.

It follows that
if $F$ is a diagonal matrix with positive eigenvalues
$s_1,\dots,s_N$,
and $\Delta F\in \M_N$ is any matrix
then the singular values $s_1',\dots,s_N'$ of $F+\Delta F$ are given by
\begin{multline} (s_i')^2=s_i^2+2 s_i \Re \Delta F_{ii} + \sum_j | \Delta F_{ji} |^2 \\ +
\sum_{j\neq i} \frac{s_i^2 |\Delta F_{ij}|^2 + 2 s_i s_j \Re (\Delta F_{ij} \Delta F_{ji})+
s_j^2 |\Delta F_{ji}|^2}{s_i^2-s_j^2}+ O\big( \| \Delta F\|^3  \big)
\label{eq:katorga2}
\end{multline}
and that the map $\Delta F\mapsto (s_1',\dots,s_N')$ is $C^2$ on some neighbourhood of $0$.

In the general case every matrix $X$ can be written as $X=U F V$, where $F$ is a
positive diagonal matrix and $U$, $V$ are unitaries.
If the singular values of $X$ are $s_1,\dots,s_N$ then (\ref{eq:katorga2}) gives us singular values
$s_1',\dots,s_N'$ of the matrix $X+\Delta X$, where $\Delta F$ is defined by
$\Delta F=V\gwia \Delta X U\gwia$.

\subsubsection{Trajectories of the Brownian motion avoid sigularities of $s$}
\label{subsubsec:avoid}
 The set of singularities of the map $s$, namely
$$\{m\in\M_N: \det m=0\mbox{ or } s_i(m)=s_j(m) \mbox{ for
some } i\neq j\},$$
is a manifold of codimension $2$ and hence almost every trajectory of a matrix Brownian
motion $A_t$ will avoid this set. In this subsection we will present a rigorous proof
of this statement.

For every $\epsilon>0$ we define a set
$$ K_\epsilon=\bigg\{ m\in \M_N: \sum_i \ln s_i(m) \geq \epsilon,
 \sum_{i<j} \ln |s_i(m)^2-s_j(m)^2|\geq \epsilon \bigg\}. $$

First of all, for any fixed $\epsilon>0$ we define a stopping time
$$T(\omega)=\min\{t\geq 0: \ln |\det A_t(\omega)|\leq \epsilon \}$$
and a stopped Brownian motion
$$\tilde{A}_t(\omega)= A_{\min[t, T(\omega)]}(\omega).$$
In the proof of Proposition \ref{prop:niespada} we showed
the function $m\mapsto \ln | \det m |$ is subharmonic and
hence $t\mapsto \ln |\det \tilde{A}_t|$ is a submartingale. It follows that
\begin{multline} \ln |\det A|  \leq
\E \ln |\det \tilde{A}_T| \leq \E \max\big(0,  \ln |\det A_T| \big)+ \\
\ln \epsilon \ P\big(\omega\in\Omega:  \ln |\det A_t(\omega)|\leq \epsilon
\mbox{ for some } 0\leq t\leq T \big) .
\label{eq:slynneszacowanie}
\end{multline}
Since the first summand on the right--hand side of the above inequality is clearly finite,
it follows that
$$\lim_{\epsilon\rightarrow 0^+} P\big(\omega\in\Omega: |\det A_t(\omega)|\geq \epsilon \mbox{ for every }
0\leq t \leq T \big) = 1.$$

Secondly, we consider a function on $M_N(\C)$ given by
\begin{equation}
m\mapsto \sum_{i<j} \ln \Big| \big(s_i(m)\big)^2 - \big(s_j(m)\big)^2 \Big| .
\label{eq:igorkraszewski}
\end{equation}
Formula (\ref{eq:katorga2}) gives us first and second derivatives of the map
$m\mapsto \big(s_1(m),\dots,s_N(m) \big)$ and allows us to find the Laplacian of
the each summand in (\ref{eq:igorkraszewski}):
$$\frac{\nabla^2 \ln |s_i^2-s_j^2|}{4} = \frac{s_i^2+s_j^2}{(s_i^2-s_j^2)^2}
+ \sum_{k\neq i,j} \frac{s_i^2+s_k^2}{(s_i^2-s_j^2)(s_i^2-s_k^2)} -
 \frac{s_j^2+s_k^2}{(s_i^2-s_j^2)(s_j^2-s_k^2)}. $$
It is not difficult to see that for every $i,j,k$ all different we have
\begin{multline*} \Bigg( \frac{s_i^2+s_k^2}{(s_i^2-s_j^2)(s_i^2-s_k^2)} -
 \frac{s_j^2+s_k^2}{(s_i^2-s_j^2)(s_j^2-s_k^2)} \Bigg) +
\Bigg( \frac{s_j^2+s_i^2}{(s_j^2-s_k^2)(s_j^2-s_i^2)} -    \\
 \frac{s_k^2+s_i^2}{(s_j^2-s_k^2)(s_k^2-s_i^2)}\Bigg) +
\Bigg( \frac{s_k^2+s_j^2}{(s_k^2-s_i^2)(s_k^2-s_j^2)} -
 \frac{s_i^2+s_j^2}{(s_k^2-s_i^2)(s_i^2-s_j^2)} \Bigg)=0 \end{multline*}
and due to these cancellations
$$\nabla^2 \sum_{i<j} \ln |s_i^2-s_j^2|= 4 \sum_{i<j} \frac{s_i^2+s_j^2}{(s_i^2-s_j^2)^2} >0 $$
holds. It follows that
$$t\mapsto  \sum_{i<j} \ln \Big| \big(s_i(\tilde{A}_t)\big)^2 - \big(s_j(\tilde{A}_t)\big)^2 \Big| $$
is a submartingale and by similar arguments as in (\ref{eq:slynneszacowanie})
we see that
$$\lim_{\epsilon\rightarrow 0^+} P\big(\omega\in\Omega: \sum_{i<j} \ln |s_i(m)^2-s_j(m)^2|\geq \epsilon
 \mbox{ for every } 0\leq t \leq T \big) = 1.$$

\subsubsection{Stochastic differential equation for $\lambda_i$}
We recall that every matrix $m$ can be written as $m=U(m) F(m) V(m)$, where $F(m)$ is a
positive diagonal matrix and $U(m)$, $V(m)$ are unitaries.
Let us define now a new matrix--valued stochastic process $B$ given by a stochastic
differential equation $dB= \big(V(A_t)\big)\gwia (dM) U(A_t)\gwia$.
It is easy to see that $B$ is again a standard matrix Brownian motion.

For any fixed $\epsilon>0$ let us consider any $C^2$ function
$\tilde{s}:\M_N \rightarrow \R^N$ such that $\tilde{s}(m)=s(m)$ for
every matrix $m\in K_\epsilon$ and such that $\|s(m)\|\leq C \|m\|$
for some universal constant $C$ and all $m\in\M_N$, where $\|\cdot \|$
denotes any norm on $\M_N$ or $\R^N$ respectively.

Function $\tilde{s}=(\tilde{s}_1,\dots,\tilde{s}_N)$ fulfills assumptions of It\^o theorem,
hence for every $T>0$ we can write the It\^o formula
\begin{multline*} \tilde{s}_i(A_T) =
 \int_0^T \sum_{k,l} \frac{\partial \tilde{s}_i(m)}{\partial m_{kl}}\bigg|_{m=A_t} dM_{kl}+   \\
\sum_{k,l} \frac{1}{4N} \left( \frac{\partial^2 \tilde{s}_i(m)}{(\partial \Re m_{kl})^2 }\bigg|_{m=A_t}+
\frac{\partial^2 \tilde{s}_i(m)}{(\partial \Im m_{kl})^2 }\bigg|_{m=A_t} \right) dt .
\end{multline*}
For every $\omega\in\Omega$ such that $A_t(\omega)\in K_\epsilon$ for all $0\leq t\leq T$ the
left--hand side of this equation is equal to $\lambda_i(T)$ and the right--hand side
can be computed from (\ref{eq:katorga2}):
$$\lambda_i(T) = \int_0^T \Re (dB_{ii})+
\frac{dt}{2 \lambda_i} \left(1-\frac{1}{2 N}+
\sum_{j\neq i} \frac{\lambda_i^2+\lambda_j^2}{N (\lambda_i^2-\lambda_j^2)}
\right). $$
Since almost every $\omega$ has the property that for some $\epsilon>0$ we have
$A_t\in K_\epsilon$ for all $0\leq t\leq T$ hence the above equation holds without
any restrictions for $\omega$. Equivalently,
$$d\lambda_i(t)=\Re (dB_{ii})+
\frac{dt}{2 \lambda_i} \left(1-\frac{1}{2 N}+
\sum_{j\neq i} \frac{\lambda_i^2+\lambda_j^2}{N (\lambda_i^2-\lambda_j^2)}
\right). $$


\subsection{Determinant of a standard Gaussian random matrix}
\begin{proposition}
\label{prop:wyznacznikgausa}
Let $(G^{(N)})$ be a sequence of independent standard Gaussian
random matrices and let $c$ be a circular element. Then
$$\lim_{N\rightarrow\infty}  \tr \ln |G^{(N)}|=
\ln \Delta(c)=- \frac{1}{2}$$
holds almost surely.

Furthermore for any $s\in\R$ we have
$$\E \min\Big(s, \inf_{N} \tr \ln |G^{(N)} | \Big) > -\infty. $$
\end{proposition}
\begin{proof}
The square of a circular element is a free Poisson element with parameter
$1$. The probability density of this element can be explicitly calculated
\cite{VDN} and the integral
$\ln \Delta(c)=\int_0^\infty \ln r \ d\mu_{\sqrt{cc\gwia}}$ can be computed directly.

Let us fix $N\in\N$. Let $v_1,\dots,v_N$ be random vectors in $\C^N$ which are
defined to be columns of the matrix $G^{(N)}$. We define
$$V_i=\sqrt{\det \left[ \begin{array}{ccc} \langle v_1, v_1 \rangle &
\cdots & \langle v_1,v_i \rangle \\
\vdots & & \vdots \\
\langle v_i,v_1 \rangle & \cdots  & \langle v_i,v_i \rangle \end{array}
\right]}, $$
 where $\langle \cdot,\cdot\rangle$ is the standard hermitian form on
$\C^N$. The above matrix $[\langle v_k,v_l\rangle]_{1\leq k,l\leq i}$ is the
complex analogue of the Gram matrix; therefore---informally
speaking---we can regard $V_i$ to be the ``complex volume'' of the
``complex parallelepiped'' defined by vectors $v_1,\dots,v_i$.

Of course $V_{i+1}$ is equal to the product of $V_i$
and $l_{i+1}$, where $l_{i+1}$ is the length of the projection of the
vector $v_{i+1}$ onto the orthogonal complement of the vectors
$v_1,\dots,v_i$. Since
$$\left[ \begin{array}{ccc} \langle v_1, v_1 \rangle &
\cdots & \langle v_1,v_N \rangle \\
\vdots & & \vdots \\
\langle v_N,v_1 \rangle & \cdots  & \langle v_N,v_N \rangle \end{array}
\right]= (G^{(N)})\gwia G^{(N)}$$
it follows that
$$| \det G^{(N)} |= V_N = l_1 l_2 \cdots l_N .$$

It is easy to see that the distribution of $l_{i}$ coincides with the
distribution of the length of a random Gaussian vector with an appropriate
covariance in the complex $(N-i+1)$--dimensional space and therefore
$$ \E l_i^{-h} = \frac{\int_0^{\infty} r^{-h} r^{2(N-i+1)-1}
e^{-N r^2}\ dr}{ \int_0^{\infty} r^{2 (N-i+1)-1} e^{-N r^2}\ dr} =
\frac{N^{\frac{h}{2}} \Gamma(N-i+1-\frac{h}{2})}{\Gamma(N-i+1)} $$
and hence Markov inequality gives us
$$P(l_N^{-1}> e^{ N \epsilon}) <  e^{-N \epsilon} \sqrt{\pi N}$$
$$ P\left[ (l_1 \cdots l_{N-1})^{-2} > e^{(1+2 \epsilon) N}\right] <
e^{(-1-2 \epsilon) N} N^{N-1} \frac{1}{\Gamma(N)}   <
e^{-2 \epsilon N}.  $$
Above we have used that random variables $l_i$ are independent and simple inequality
$\Gamma(N)> \left( \frac{N-1}{e}  \right)^{N-1} $ for $N\in\N$.

Since
\begin{multline*} P\left( \frac{\log l_1+\cdots+\log l_N}{N} < -\frac{1}{2} - 2 \epsilon \right) \leq   \\
 P\left( \frac{\log l_1+\cdots+\log l_{N-1}}{N} < -\frac{1}{2}-\epsilon \right) +
 P\left( \frac{\log l_N}{N} < -\epsilon\right) \end{multline*}
Borel--Cantelli lemma implies
$$\liminf_{N\rightarrow\infty} \tr \ln |G^{(N)}| \geq \ln \Delta(c)$$
almost surely. This together with Lemma \ref{lem:kopalny} gives
us the first part of the proposition.

It is possible to find a constant $C$ such that for every $\epsilon>\frac{1}{4}$ we have
$$P\left(\inf_N \tr \ln |G^{(N)}| < -\frac{1}{2} -2\epsilon \right) \leq
\sum_N  P\left(\tr \ln |G^{(N)}| < -\frac{1}{2} -2\epsilon \right) \leq C e^{-\epsilon}.$$
If $\nu$ is the distribution of the random variable $\inf_N \tr \ln |G^{(N)}|$ then
integration by parts gives
$$\int_{-\infty}^{-1} (t+1) d\nu(t)= - \int_{-\infty}^{-1} \nu(-\infty,t)\ dt > -\infty$$
and the second part follows.
\end{proof}


\section{Acknowledgements}
I would like to thank Roland Speicher for many fruitful discussions.
Also, I would like to thank the Reviewer for a very careful reading of the article and
numerous very helpful remarks. A part of the research was conducted at Texas A\&{}M University on a
scholarship funded by Polish--US Fulbright Commission.
I acknowledge the support of Polish Research Committee
grant No.\ P03A05415.

\end{document}